constants). Returning to the original variables, we obtain a two-parameter family of two-frequency quasiperiodic solutions of (22):

$$x(t) = c_1 \cos t + c_2 \sin t, \quad z(t) = 0,$$

$$y(t) = -\frac{1}{2}(c_1^2 + c_2^2) - c_1 \sin t + c_2 \cos t - \alpha c_1 c_2 \sin 2t - \frac{1}{2}\alpha(c_1^2 - c_2^2)\cos 2t,$$

the basis frequency of these solutions in conjunction with the frequencies of the right side form a rationally linearly independent set.

## LITERATURE CITED

1. É. I. Grudo, Differents. Uravn., 22, No. 9, 1499-1504 (1986).
2. A. K. Demenchuk, Differents. Uravn., 26, No. 12, 2167-2170 (1990).
3. Yu. N. Ribikov, Differents. Uravn., 7, No. 8, 1347-1356 (1971).

# INTEGRAL ESTIMATES FOR ORDINARY DIFFERENTIAL EQUATIONS AND THEIR APPLICATION TO NONSMOOTH OPTIMAL-CONTROL PROBLEMS

N. G. Dokuchaev                                                                                     UDC 517.977

Introduction. We investigate the differential equation

$$\frac{dy}{dt}(t) = f[y(t), t]. \tag{0.1}$$

and write $y^{x,s}(t)$ for a solution of (0.1) on $[s, T]$ with the initial condition

$$y(s) = x. \tag{0.2}$$

Here $0 \leq s < t \leq T$ with $T > 0$ a given number, and $x \in \mathbf{R}^n$. The function $f(x, t): \mathbf{R}^n \times [0, T] \to \mathbf{R}^n$ and its x-derivatives up to the second order inclusive are measurable (the derivatives exist for each t).

We investigate functionals

$$V(x, s) = \int_s^T \varphi[y^{x,s}(t), t] dt, \tag{0.3}$$

of solutions of (0.1), (0.2), where $\varphi(x, t): \mathbf{R}^n \times [0, T] \to \mathbf{R}$ is Lebesgue measurable.

If $\varphi(x, t)$ is x-differentiable, then V(x, s) coincides with the solution of the corresponding Cauchy problem for a first-order partial differential equation [problem (1.1)]; this solution of the boundary-value problem will be denoted by v(x, s). It is known ([1], p. 649) that

$$\int_{\mathbf{R}^n} |v(x, s)|^2 dx \leq c \int_s^T dt \int_{\mathbf{R}^n} |\varphi(x, t)|^2 dx \tag{0.4}$$

with the same constant c > 0 for all $s \in [0, T]$ and $\varphi \in C^1(\mathbf{R}^n \times [0, T]) \cap L_2(\mathbf{R}^n \times [0, T])$. Let L be the operator mapping $\varphi$ into the solution of the boundary-value problem: $v = L\varphi$. By virtue of (4) this operator can be extended by continuity to a continuous operator $L: L_2(\mathbf{R}^n \times [0, T]) \to C([0, T] \to L_2(\mathbf{R}^n))$. However it is not obvious whether $v = L\varphi$ with V as in (0.3) for arbitrary $\varphi \in L_2(\mathbf{R}^n \times [0, T])$. To prove this is similar to proving, for such $\varphi$, that (0.4) holds for the functional (0.3). We continue this type of investigation. We prove, for an arbitrary function $\varphi$ from the class $L^1([0, T] \to L_2(\mathbf{R}^n))$, that V(x, s) is in $C([0, T] \to L_2(\mathbf{R}^n))$ and that V = v where $v = L\varphi$, and derive the inequality

$$\left(\int_{\mathbf{R}^n} |V(x, s)|^2 dx\right)^{1/2} \leq c \int_s^T dt \left(\int_{\mathbf{R}^n} |\varphi(x, t)|^2 dx\right)^{1/2} \tag{0.5}$$

with a constant $c=c(s)=\exp\int_s^T \delta_f(t)\,dt$, where $\delta_f(t)=\frac{1}{2}\sup_{x\in\mathbf{R}^n}\left|\sum_{i=1}^n \frac{\partial f_i}{\partial x_i}(x,t)\right|$, and $x_i$ and $f_i$ are the components of x and f.

We use (0.5) to investigate a control problem for an ordinary differential equation for the minimization of a functional (0.2) with a nondifferentiable $\varphi$. There are many articles on such problems (see the review in [2]). We first avoid the case of an arbitrary Lebesgue measurable function $\varphi$ (which cannot have a quasiderivative or a subdifferential, being discontinuous and unbounded). There are several methods of "smoothing" the problem: the original equation of the problem can be replaced by a nonsingular Ito equation (which lowers the smoothness required of $\varphi$), or we can approximate $\varphi$ by smooth functions. Here we use another approach: we show that the mean value of the functional (0.3) will depend "smoothly" on the variations of the control, if the initial state vector is assumed to be random with density from the class $W_2^1(\mathbf{R}^n)$. This method of "smoothing" a nonsmooth problem essentially conserves the equation of the problem and the quality functional. For this formulation of the problem we find necessary conditions for optimization, and sufficient conditions for the existence of an optimal control.

In the following all vectors and matrices are real, and M indicates mathematical expectation.

1. Integral Estimates of Trajectories. We introduce the notation for the Banach spaces used. For integers $k \geq 0$ we write $H^k$ for the standard Hilbert spaces $H^k = W_2^k(\mathbf{R}^n)$. We also introduce the Banach spaces $\mathscr{C}^k = C([0, T] \to H^k)$ (of continuous functions of t with values in $H^k$), $X^k = L^2([0, T] \to H^k)$, and $X_1^k = L^1([0, T] \to H^k)$. Elements of these spaces are Lebesgue measurable functions $u(x, t): \mathbf{R}^n \times [0, T] \to \mathbf{R}$, satisfying the appropriate smoothness and summability conditions (although, strictly speaking, the elements of $X^0$ are not functions but equivalence classes). For a space $\mathscr{X}$ of functions $u(x, t): \mathbf{R}^n \times [0, T] \to \mathbf{R}$, we write $\mathscr{X}[s, T]$ for the space of restrictions of functions from $\mathscr{X}$ to $\mathbf{R}^n \times [s, T]$. For the space $\mathscr{X}$, $\|\cdot\|_{\mathscr{X}}$ denotes the norm; for a Hilbert space $(\cdot, \cdot)_{\mathscr{X}}$ denotes the scalar product.

It is known that, if $\varphi \in C([0, T] \to C^1(\mathbf{R}^n)) \cap X^1$, then the integral (0.3) of V(x, s) coincides with the solution v(x, s) of the Cauchy problem

$$\frac{\partial v}{\partial s}(x, s) + A(x, s)v(x, s) = -\varphi(x, s), \quad v(x, T) = 0. \tag{1.1}$$

Here $0 \leq s \leq T$, and the differential operator A(x, t) is as in the relation

$$A(x, t)v(x) = \sum_{i=1}^n f_i(x, t)\frac{\partial v}{\partial x_i}(x) \tag{1.2}$$

where $x_i$ and $f_i$ are the components of x and f.

Let $A^*(x, t)$ be the conjugate of A(x, t):

$$A^*(x, t)v(x) = -\sum_{i=1}^n \frac{\partial}{\partial x_i}(f_i(x, t)v(x)).$$

To prove that V = v, we consider problem (1.1) in conjugation with the Cauchy problem

$$\frac{\partial p}{\partial t}(x, t) = A^*(x, t)p(x, t), \quad p(x, s) = \rho_s(x), \tag{1.3}$$

where $0 \leq s < t \leq T$ and $\rho_s \in H^1$.

Let L and $\mathscr{L}_s^*$ map $\varphi$ and $\rho_s$ into the solutions $v = L\varphi$ and $p = \mathscr{L}_s^*\rho_s$ of problems (1.1) and (1.3), respectively.

It is known that the operators $L: X^1 \to \mathscr{C}^0$, $L: X^1 \to X^1$, $\mathscr{L}_s^*: H^1 \to X^1[s, T]$, and $\mathscr{L}_s^*: H^1 \to \mathscr{C}^0[s, T]$ are continuous. This is a special case of Theorem 4.2.1 from [3] concerning more general equations.

Let $a = a(\omega)$ be a random vector (here $\omega$ is an abstract variable characterizing the case; we assume that the probability space $(\Omega, \mathscr{F}, P)$ is given ([3], p. 15), and write $\omega \in \Omega$ for an elementary event). Let $y^{a,s}(t)$ denote the random process satisfying (0.1) for each $\omega$ and the initial condition

$$y(s) = a. \tag{1.4}$$

We assume that the vector a has distribution density $\rho_s(X)$. If $\rho_s \in H^1$, then it is known that $p(x, t) = \mathscr{L}_s^*\rho_s \in X^1[s, T] \cap \mathscr{C}^0[s, T]$ is the density distribution of the process $y^{a,s}(t)$ (this follows from Theorem 5.3.1 [3] for more general processes). In the sequel we extend this result to the somewhat more complicated case $\rho_s \in H^0$, in which it is formally impossible to use Eq. (1.3) because p does not have derivatives.

The following lemma contains a specific definition needed in our problem, in terms of known results.

LEMMA 1.1. The operators L and $\mathscr{L}_s^*$ can be continued to continuous linear operators $L: X_1^0 \to \mathscr{C}^0$ and $\mathscr{L}_s^*: H^0 \to \mathscr{C}^0[s, T]$ with the norm of the extended L not exceeding $\exp \int_0^T \delta_f(t)\,dt$ and the norm of the extended $\mathscr{L}_s^*$ not exceeding $\exp \int_s^T \delta_f(t)\,dt$, where

$$\delta_f(t) = \frac{1}{2} \sup_{x \in \mathbf{R}^n} \left| \sum_{i=1}^n \frac{\partial f_i}{\partial x_i}(x, t) \right|.$$

In the following L and $\mathscr{L}_s^*$ are understood to be the operators in Lemma 1.1 if nothing is said to the contrary.

LEMMA 1.2. If $s \in [0, T]$, $\varphi \in X_1^0$, $\rho_s \in H^0$, $v = L\varphi \in \mathscr{C}^0$, and $p = \mathscr{L}_s^* \rho_s \in \mathscr{C}^0[s, T]$, then

$$(v(x, s), \rho_s(x))_{H^0} = \int_s^T (\varphi(x, t), p(x, t))_{H^0}\, dt. \tag{1.5}$$

THEOREM 1.1. If $\varphi \in X_1^0$ is Lebesgue measurable, then $V(x, s)$, defined by (0.3), is Lebesgue measurable with respect to x, s and, for each $s \in [0, T]$, is uniquely determined as an element of $L_2(\mathbf{R}^n)$, i.e., for almost all $x \in \mathbf{R}^n$. Moreover $V \in \mathscr{C}^0$ and $V = v$ in $\mathscr{C}^0$, where $v = L\varphi \in \mathscr{C}^0$ with the operator L as in Lemma 1.1.

THEOREM 1.2. If $a = a(\omega)$ is a random n-vector with density distribution $\rho_s(x) \in H^0$, then $p(x, t) = \mathscr{L}_s^* \rho_s \in \mathscr{C}^0[s, T]$ is the distribution density of the process $y^{a,s}(t)$ determined by Eq. (0.1) and the initial condition (1.4).

Proof of Lemma 1.1. We apply the customary method of choosing the energy inequalities [4] for the norm in $H^0 = L_2(\mathbf{R}^n)$ of solutions v and p of the Cauchy problems (1.1) and (1.3).

For $u \in H^1$, the following relations hold for almost all t:

$$(u(x), A(x, t)u(x) + \varphi(x, t))_{H^0} = \sum_{i=1}^n \left( u(x), f_i(x, t) \frac{\partial u}{\partial x_i}(x) \right)_{H^0} +$$

$$+ (u(x), \varphi(x, t))_{H^0} = \frac{1}{2} \sum_{i=1}^n \left( \frac{\partial}{\partial x_i}[u(x)^2], f_i(x, t) \right)_{H^0} +$$

$$+ (u(x), \varphi(x, t))_{H^0} = -\frac{1}{2} \left( u(x)^2, \sum_{i=1}^n \frac{\partial f_i}{\partial x_i}(x, t) \right)_{H^0} +$$

$$+ (u(x), \varphi(x, t))_{H^0} \leq \delta_f(t) \|u(x)\|_{H^0}^2 + \|\varphi(x, t)\|_{H^0} \|u(x)\|_{H^0}.$$

Let $\varphi \in X^1$; then $v = L\varphi \in \mathscr{C}^0 \cap X^1$. Reasoning as in [4] (p. 160), we conclude that $\|v(x, t)\|_{H^0}^2$ is an absolutely continuous function of t:

$$\|v(x, t)\|_{H^0}^2 = 2 \int_t^T (v(x, s), A(x, s)v(x, s) + \varphi(x, s))_{H^0}\, dt.$$

In other words,

$$\frac{\partial}{\partial t} \|v(x, t)\|_{H^0}^2 = 2(v(x, s), A(x, s)v(x, s) + \varphi(x, s))_{H^0}$$

for almost all t; hence

$$\frac{\partial}{\partial t} (\|v(x, t)\|_{H^0}^2) \leq 2[\delta_f(t) \|v(x, t)\|_{H^0} + \|\varphi(x, t)\|_{H^0}] \|v(x, t)\|_{H^0}$$

for almost all t. Hence, for almost all t, we have

$$\frac{\partial}{\partial t} (\|v(x, t)\|_{H^0}^2) = 2\|v(x, t)\|_{H^0} \frac{\partial}{\partial t} \|v(x, t)\|_{H^0}.$$

Hence $\|v(x, t)\|_{H^0}$ is an absolutely continuous function of t and, for almost all t,

$$\frac{\partial}{\partial t}\|v(x, t)\|_{H^0} \leqslant \delta_f(t)\|v(x, t)\|_{H^0} + \|\varphi(x, t)\|_{H^0}.$$

Results in [4] (p. 152) now imply that

$$\|v(x, s)\|_{H^0} \leqslant \left(\exp \int_s^T \delta_f(t)\, dt\right)\|\varphi\|_{X_1^0[s, T]}.$$

The conclusion of the lemma for L is obtained by extending this inequality to $\varphi \in X_1^0$, with a set $X^1$ everywhere dense in $X_1^0$.

The conclusions of the lemma for the operators $\mathcal{L}_s^*$ are proved similarly. The inequality $(p(x), A^*(x, t)p(x))_{H^0} \leqslant \delta_f(t)\|p(x)\|_{H^0}$ ($\forall p \in H^1$) and the relation $\|p(x, t)\|_{H^0}^2 = 2\int_s^t (p(x, \tau), A^*(x, \tau)p(x, \tau))_{H^0} d\tau + \|p(x, s)\|_{H^0}$ are used. This proves Lemma 1.1.

<u>Proof of Lemma 1.2.</u> First assume that $\varphi \in X^1$ and $\rho_s \in H^1$; then Theorem 4.2.1 [3] implies that $\partial v(x, t)/\partial t$ and $A(x, t)v(x, t)$ are in $X^0$, $\partial p(x, t)/\partial t$ and $A^*(x, t)p(x, t)$ are in $X^0[s, T]$, and $v \in X^1$ and $p \in X^1[s, T]$. In this case (1.5) follows from the relations

$$-(v(x, s), \rho_s(x))_{H^0} = (v(x, T), p(x, T))_{H^0} - (v(x, s), \rho_s(x))_{H^0} =$$
$$= \int_s^T \frac{d}{dt}(v(x, t), p(x, t))_{H^0} dt = \int_s^T \{(-A(x, t)v(x, t) -$$
$$-\varphi(x, t), p(x, t))_{H^0} + (v(x, t), A^*(x, t)p(x, t))_{H^0}\} dt = -\int_s^T (\varphi(x, t), p(x, t))_{H^0} dt.$$

For arbitrary $\varphi \in X_1^0$ and $\rho_s \in H^0$, we construct approximating sequences $\varphi^{(i)} \in X^1$ and $\rho_s^{(i)} \in H^1$ $i = 1, 2, \ldots$, such that $\varphi^{(i)} \to \varphi$ in $X_1^0$ and $\rho_s^{(i)} \to \rho_s$ in $H^0$ for $i \to +\infty$. From Lemma 1.1 we have $v^{(i)} = L\varphi^{(i)} \to v$ in the metric of $\mathscr{C}^0$ and $p^{(i)} = \mathcal{L}_s^* \rho_s^{(i)} \to p$ in the metric of $\mathscr{C}^0[s, T]$ for $i \to +\infty$. The conclusion of Lemma 1.2 follows.

<u>Proof of Theorem 1.1.</u> Let $\rho_s \in H^1$ be the distribution density of the random n-vector $a = a(\omega)$. With no loss of generality we assume that $a = a(\omega)$ is a random vector on the probability space $(\Omega, \mathscr{F}, P)$, where $\Omega = \{\omega\} = \mathbf{R}^n$, and $\mathscr{F}$ is the completion of the σ-algebra of Borel sets in $\mathbf{R}^n$, and P is a probability measure on $\mathscr{F}$ such that $P(d\omega) = \rho_s(\omega)d\omega$; we also assume that $a(\omega) \equiv \omega$.

We prove that, for each s, $\varphi[y^{a(\omega),s}(t), t]: \mathbf{R}^n \times [s, T] \to \mathbf{R}$ is Lebesgue measurable as a function of $\omega$, t. If $\varphi(x, t)$ is Borel measurable, then this follows directly from properties of a composition of measurable functions ([5], p. 282) and the Lebesgue measurability of a mapping $Q(\omega, t) = [y^{a(\omega),s}(t), t]: \mathbf{R}^n \times [s, T] \to \mathbf{R}^{n+1}$ [which follows from the continuous dependence of the solution of the Cauchy problem (0.1), (0.2) on x, t].

Now let $\varphi(x, t)$ be only Lebesgue measurable, let $\lambda_m(\cdot)$ be the Lebesgue measure in $\mathbf{R}^m$, and let $\mathscr{B}_m$ and $\bar{\mathscr{B}}_m$ be the sets of Borel and Lebesgue measurable sets in $\mathbf{R}^m$, respectively. We use the notation $\mathscr{B}_m^0 = \{e \in \mathscr{B}_m : \lambda_m(e) = 0\}$ and $\bar{\mathscr{B}}_m^0 = \{e \in \bar{\mathscr{B}}_m : \lambda_m(e) = 0\}$. We note that, for each t and $e \in \mathscr{B}_n^0$, the inclusion $Q^{-1}(\cdot, t)(e) \in \mathscr{B}_n^0$ holds [this is a phase-volume conservation property of Eq. (0.1)]. Let $N_L \in \bar{\mathscr{B}}_{n+1}^0$ be an arbitrary set and let $N_L \subset \mathbf{R}^n \times [s, T]$. Then there is $N_B \in \mathscr{B}_{n+1}$ such that $N_L \subset N_B \subset \mathbf{R}^n \times [s, T]$. Using the notation $N_B(t) = \{x \in \mathbf{R}^n : (x, t) \in N_B\}$, we find that $\lambda_n(N_B(t)) = 0$ and $Q^{-1}(\cdot, t)(N_B(t)) \in \mathscr{B}_n^0$ for almost all t. Moreover $Q^{-1}(N_B) \in \mathscr{B}_{n+1}$, and so $\lambda_{n+1} \times (Q^{-1}(N_B)) = \int_s^T \lambda_n(Q^{-1}(\cdot, t)[N_B(t)])\, dt = 0$. Since $Q^{-1}(N_L) \subset Q^{-1}(N_B)$ we conclude that $Q^{-1}(N_L) \in \bar{\mathscr{B}}_{n+1}^0$ for arbitrary $N_L \in \bar{\mathscr{B}}_{n+1}^0$, $N_L \subset \mathbf{R}^n \times [s, T]$.

Hence $Q^{-1}(\Gamma) \in \bar{\mathscr{B}}_{n+1}$ for arbitrary $\Gamma \in \bar{\mathscr{B}}_{n+1}$, $\Gamma \subset \mathbf{R}^n \times [s, T]$ [in fact there is a representation $\Gamma = \Gamma_B \cup N_L$ where $\Gamma_B \in \mathscr{B}_{n+1}$ and $N_L \in \bar{\mathscr{B}}_{n+1}^0$, and for which $Q^{-1}(\Gamma) = Q^{-1}(\Gamma_B) \cup Q^{-1}(N_L)$, and where $Q^{-1}(\Gamma_B) \in \mathscr{B}_{n+1}$ and $Q^{-1}(N_L) \in \bar{\mathscr{B}}_{n+1}$]. The usual method now shows that the composition $\varphi[Q(\omega, t)] = \varphi[y^{a(\omega),s}(t), t]$ is Lebesgue measurable as a function of $\omega$, t ([5], p. 282).

We have already shown that $p(x, t) = \mathcal{L}_s^* \rho_s \in \mathscr{C}^0[s, T] \cap X^1[s, T]$ is the density distribution of a random process $y^{a,s}(t)$ determined by (0.1) and (0.4). Hence, for almost all t,

$$(\varphi(x, t), p(x, t))_{H^0} = M\varphi[y^{a,s}(t), t] \tag{1.6}$$

(here we can use Proposition 3 ([6], p. 288) and Theorem 1 ([6], p. 331)). By virtue of Lemma 1.2 for $v = L\varphi \in \mathscr{C}^0$, we conclude from (1.5) that

$$(v(x,s), \rho_s(x))_{H^0} = \int_s^T (\varphi(x,t), p(x,t))_{H^0} dt = \int_s^T M\varphi[y^{a,s}(t), t] dt. \qquad (1.7)$$

We note further that $M|\varphi[y^{a(\omega),s}(t), t]| \in L_1(s, T)$, since the function $|\varphi(x,t)|$ can be substituted in (1.6), and so $M|\varphi[y^{a(\omega),s}(t), t]| = (|\varphi(x,t)|, p(x,t))_{H^0} \in L_1(s, T)$.

Fubini's theorem can be applied, and we therefore obtain

$$\int_s^T M\varphi[y^{a(\omega),s}(t), t] dt = M \int_s^T \varphi[y^{a(\omega),s}(t), t] dt = \int_\Omega P(d\omega) \int_s^T \varphi[y^{a(\omega),s}(t), t] dt = (V(x,s), \rho_s(x))_{H^0}$$

by virtue of the choice of the measure $P(d\omega)$. Hence (1.7) implies that

$$(v(x,s), \rho_s(x))_{H^0} = (V(x,s), \rho_s(x))_{H^0} \qquad (1.8)$$

for each $s \in [0, T]$ and arbitrary $\rho_s \in H^1$ such that: a) $\rho_s(x) \geq 0$ for almost all x; b) $\rho_s \in L_1(\mathbf{R}^n)$; c) $\int_{\mathbf{R}^n} \rho_s(x) dx = 1$. The usual reasoning now establishes (1.8) for $\rho_s \in H^1$. Thus $V(x,s) = v(x,s)$ in $H^0$ for each s. Lemma 1.1 implies that $V \in \mathscr{C}^0$ and $v = V$ in $\mathscr{C}^0$ and this proves Theorem 1.1.

Proof of Theorem 1.2. Let $\tilde{\varphi}(x) \in C(\mathbf{R}^n) \cap H^1$ and $\xi(t) \in L_2(s, T)$ be arbitrary. For $\varphi(x,t) = \tilde{\varphi}(x)\xi(t)$ and $v = L\varphi$ the following relations hold:

$$Mv[a(\omega), s] = (v(x,s), \rho_s(x))_{H^0} = M\int_s^T \tilde{\varphi}[y^{a(\omega),s}(t)]\xi(t) dt = \int_s^T M\tilde{\varphi}[y^{a(\omega),s}(t)]\xi(t) dt.$$

The applicability of Fubini's theorem in the third relation here is proved as in the proof of Theorem 1.1. It follows from (1.5) that

$$\int_s^T M\tilde{\varphi}[y^{a(\omega),s}(t)]\xi(t) dt = \int_s^T \xi(t)(\tilde{\varphi}(x), p(x,t))_{H^0} dt$$

and since $\xi$ is arbitrary, we conclude that

$$M\tilde{\varphi}[y^{a(\omega),s}(t)] = (\tilde{\varphi}(x), p(x,t))_{H^0} \qquad (1.9)$$

for almost all $t \in [s, T]$. We know that $p(x, t)$ is a continuous function of t in $H^0$ and the process $y^{a,s}(t)$ is continuous in $L^2(\Omega, \mathscr{F}, P, \mathbf{R}^n)$ (this is easily proved). Hence (1.9) holds for for all $t \in [s, T]$, and since $\tilde{\varphi}$ is arbitrary the conclusion of Theorem 1.2 holds.

2. Nonsmooth Optimal-Control Problem. We consider the following optimal-control problem:

$$\frac{dy}{dt}(t) = f[y(t), u(t), t], \qquad (2.1)$$

$$y(0) = a, \qquad (2.2)$$

$$\Phi[u(\cdot)] = M\int_0^T \varphi[y(t), u(t), t] dt = \min_{u(\cdot)}. \qquad (2.3)$$

Here $y(t)$ is the n-dimensional vector-valued state; a is the given initial state, in general random; $u(t)$ is the m-dimensional control vector.

An admissible control is a Lebesgue measurable function $u(t): [0, T] \to \mathbf{R}^m$ with values in a given set $\Delta \subset \mathbf{R}^m$ for almost all t. We write $\mathscr{U}_0$ for the set of admissible controls. These are "programmed" controls, i.e., controls without feedback, and nonrandom functions even when a and $y(t)$ are random.

The following conditions are assumed to be satisfied.

I. The function $f = f(x, u, t): \mathbf{R}^n \times \mathbf{R}^m \times [0, T] \to \mathbf{R}^n$ is measurable, satisfies a Lipschitz condition with respect to $x \in \mathbf{R}^n$ and $u \in \Delta$ uniformly with respect to t, and, together with its partial x-derivatives up to the second order inclusive (which exist for each u, t), is bounded uniformly with respect to x, u, and t. Moreover $f(x, u, t)$ and $\partial f(x, u, t)/\partial x$ are continuous functions of t for each x and u.

II. The function $\varphi = \varphi(x, u, t): \mathbf{R}^n \times \mathbf{R}^m \times [0, T] \to \mathbf{R}$ is Lebesgue measurable. Thus there exists a constant $K_0 > 0$, such that

$$\int_0^T dt \left( \int_{\mathbf{R}^n} |\varphi(x, u(t), t)|^2 dx \right)^{1/2} \leq K_0$$

for all $u(\cdot) \in \mathscr{U}_0$.

We note that the last inequality holds, for example, if $|\varphi(x,u,t)| \leq \varphi_*(x,t)$ $(\forall u \in \Delta)$, where $\varphi_* \in X_1^0$ (here and below we use the notation for spaces introduced in part 1).

For $x \in \mathbf{R}^n$, $u \in \Delta$, and $t \in [0,T]$, we introduce the differential operators defined as follows:

$$A(x,u,t)v(x) = \sum_{i=1}^n f_i(x,u,t) \frac{\partial v}{\partial x_i}(x),$$

$$A^*(x,u,t)v(x) = -\sum_{i=1}^n \frac{\partial}{\partial x_i}[f_i(x,u,t)v(x)],$$

(2.4)

where $x_i$ and $f_i$ are the components of the vectors x and f.

We also introduce the following operators defined for each $u(\cdot) \in \mathcal{U}_0$:

$$L[u(\cdot)]: X_1^0 \to \mathscr{C}^0, \quad \mathscr{L}_s^*[u(\cdot)]: H^1 \to \mathscr{C}^0[s,T], \quad \mathscr{L}_s[u(\cdot)]: H^1 \to X^1[s,T];$$

these operators are defined by the method described in Part 1 for $A(x,t) = A(x, u(t), t)$ and $A^*(x,t) = A^*(x, u(t), t)$ [see (1.2) with (2.4)]: they coincide with the operators $L: X_1^0 \to \mathscr{C}^0$, $\mathscr{L}_s^*: H^1 \to \mathscr{C}^0[s,T]$, and $\mathscr{L}_s: H^1 \to X^1[s,T]$ for $A(x,t) = A(x, u(t), t)$ and $A^*(x,t) = A^*(x, u(t), t)$.

<u>THEOREM 2.1.</u> Let a be a random vector in $\mathbf{R}^n$ with distribution density $\rho_0(x) \in W_2^1(\mathbf{R}^n)$, and let $\hat{y}(\cdot)$ and $\hat{u}(\cdot)$ be optimal for problem (2.1)-(2.3). Write $\hat{L} = L[\hat{u}(\cdot)]$ and $\hat{\mathscr{L}}_s^* = \mathscr{L}_s^*[\hat{u}(\cdot)]$ for the operators corresponding to $u(\cdot) = \hat{u}(\cdot)$, let $\hat{\varphi}(x,t) = \varphi(x, \hat{u}(t), t)$ be the corresponding function, and let $\hat{y}^{x,s}(t)$ be the solution of Eq. (2.1) with $u(\cdot) = \hat{u}(\cdot)$ and the Cauchy condition $y(s) = x$, where $x \in \mathbf{R}^n$. Then $\hat{y}(t)$ has density distribution $\hat{p}(x,t) = \hat{\mathscr{L}}_0 \rho_0 \in X^1 \cap \mathscr{C}^0$, and $\hat{v}(x,s) = \int_s^T \varphi[\hat{y}^{x,s}(t), \hat{u}(t), t] dt$ is such that $\hat{v} = \hat{L}\hat{\varphi} \in \mathscr{C}^0$ (see Theorem 1.1). The following inequality holds for each $\bar{u} \in \Delta$ for these $\hat{p}$ and $\hat{v}$:

$$\int_{\mathbf{R}^n} \left\{ \hat{v}(x,t) \sum_{i=1}^n \frac{\partial}{\partial x_i}[\hat{p}(x,t) f_i(x, \hat{u}(t), t)] - \varphi(x, \hat{u}(t), t) \hat{p}(x,t) \right\} dx \geq$$

$$\geq \int_{\mathbf{R}^n} \left\{ \hat{v}(x,t) \sum_{i=1}^n \frac{\partial}{\partial x_i}[\hat{p}(x,t) f_i(x, \bar{u}, t)] - \varphi(x, \bar{u}, t) \hat{p}(x,t) \right\} dx \quad (2.5)$$

for almost all $t \in [0,T]$; more precisely, for $t \in \Theta_1 \cap \Theta_2(\bar{u})$, where $\Theta_1 \subset [0,T]$ is the set of Lebesgue points of the function $(\hat{v}(x,t), A^*(x, \hat{u}(t), t)\hat{p}(x,t))_{H^0} + (\varphi(x, \hat{u}(t), t), \hat{p}(x,t))_{H^0}$, and $\Theta_2(\bar{u})$ is the set of Lebesgue points of $(\varphi(x, \bar{u}, t), \hat{p}(x,t))_{H^0}$. If

$$\mathrm{mes}\left([0,T] \setminus \bigcap_{\bar{u} \in \Delta} \Theta_2(\bar{u})\right) = 0, \quad (2.6)$$

then condition (2.5) is satisfied for almost all $t \in [0,T]$ and arbitrary $\bar{u} \in \Delta$.

We note that condition (2.6) is always satisfied for functions $\varphi(x, \bar{u}, t) \equiv \varphi(x,t)$ not depending on $\bar{u}$, or which $\varphi(x, \bar{u}, t) \in \mathscr{C}^0$ $(\forall \bar{u} \in \Delta)$.

It should be noted that, if $\varphi$ and $\hat{v}$ are smooth, then integration by parts in condition (2.5) yields

$$M\{\psi(t)^\mathsf{T} f(y(t), \hat{u}(t), t) - \varphi(y(t), \hat{u}(t), t)\} \geq M\{\psi(t)^\mathsf{T} f(y(t), \bar{u}, t) - \varphi(y(t), \bar{u}, t)\},$$

where $\psi(t) = -\frac{\partial \hat{v}}{\partial x}(\hat{y}(t), t)$, i.e., it becomes similar to Pontryagin's maximum principle for an optimal-control problem on the ensemble of trajectories of the type (2.1)-(2.3) (see [7]).

<u>Proof of Theorem 2.1.</u> Let $\mu = \|\bar{t}, \bar{u}\|$, where $\bar{t} \in (0,T)$ and $\bar{u} \in \Delta$. For each $\mu$ consider the interval $\Theta(\varepsilon) = [\bar{t}, \bar{t}+\varepsilon]$ for $\varepsilon \geq 0$, and let $\varepsilon_\mu > 0$ be small enough to ensure that $\Theta(\varepsilon) \subset [0,T]$ when $\varepsilon \in [0, \varepsilon_\mu]$.

For each $\mu$ and for $\varepsilon \in [0, \varepsilon_\mu]$, specify the variation $u(\cdot, \varepsilon|\mu)$ of $\hat{u}(\cdot)$, i.e., specify mappings $u(\cdot|\mu): [0, \varepsilon_\mu] \to \mathcal{U}_0$. It is sufficient to use a bundle of simple needle-shaped variation [8]: let $u(\cdot, 0|\mu) = \hat{u}(\cdot)$; if $\varepsilon > 0$ then put $u(t, \varepsilon|\mu) = \hat{u}(t)$ for $t \in [0,T] \setminus \Theta(\varepsilon)$ and $u(t, \varepsilon|\mu) = \bar{u}$ for $t \in \Theta(\varepsilon)$. We thus have a bundle of curves in the space $\mathcal{U}_0$, or a bundle of variations. To prove the theorem we find the differential of the functional (2.3) along a curve $u(\cdot, \varepsilon|\mu)$ [8].

We fix $\mu = \|\bar{t}, \bar{u}\|$, and use the following notation: $\widehat{A}(x, t) = A(x, \hat{u}(t), t)$, $A_\varepsilon(x, t) = A(x, u(t, \varepsilon|\mu), t)$, $\varphi_\varepsilon(x, t) = \varphi(x, u(t, \varepsilon|\mu), t)$, $v_\varepsilon = L[u(\cdot, \varepsilon|\mu)]\varphi_\varepsilon$, and $p_\varepsilon = \mathscr{L}_0^*[u(\cdot, \varepsilon|\mu)]\rho_0$. For $u(\cdot) = u(\cdot, \varepsilon|\mu)$ and $u(\cdot) = \hat{u}(\cdot)$ we denote the solution of problem (2.1), (2.2) by $y_\varepsilon(t)$ and $\hat{y}(t)$. The arguments $x$ or $t$ of functions from $X_1^0$ will be sometimes omitted for brevity.

By virtue of (2.3), the application of Fubini's theorem, justified as in the proof of Theorem 1.1 yields

$$\Phi[u(\cdot, \varepsilon|\mu)] - \Phi[\hat{u}(\cdot)] = \int_0^T M\varphi[y_\varepsilon(t), t]\, dt - \int_0^T M\varphi[\hat{y}(t), t]\, dt = (v_\varepsilon(x, 0), \rho_0(x))_{H^0} - (\hat{v}(x, 0), \rho_0(x))_{H^0}. \quad (2.7)$$

Here the second relation follows from Theorem 1.2 and Lemma 1.1.

Suppose that $\hat{\varphi} \in X^1$ and $\varphi_\varepsilon \in X^1$; then $v_\varepsilon \in X^1$, $\hat{v} \in X^1$, $\partial v/\partial t \in X^0$, and $\partial \hat{v}/\partial t \in X^0$ and (2.7) can be written as

$$\Phi[u(\cdot, \varepsilon|\mu)] - \Phi[\hat{u}(\cdot)] = (\hat{v}(x, T), \hat{p}(x, T))_{H^0} - (\hat{v}(x, 0), \hat{p}(x, 0))_{H^0} -$$

$$- (v_\varepsilon(x, T), \hat{p}(x, T))_{H^0} + (v_\varepsilon(x, 0), \hat{p}(x, 0))_{H^0} = \int_0^T \left\{ \left(\hat{v}, \frac{\partial \hat{p}}{\partial t}\right)_{H^0} + \right.$$

$$\left. + \left(\frac{\partial \hat{v}}{\partial t}, \hat{p}\right)_{H^0} - \left(v_\varepsilon, \frac{\partial \hat{p}}{\partial t}\right)_{H^0} - \left(\frac{\partial v_\varepsilon}{\partial t}, \hat{p}\right)_{H^0} \right\} dt =$$

$$= \int_0^T \{(\hat{v}, \widehat{A}^*\hat{p})_{H^0} + (-\widehat{A}\hat{v} - \hat{\varphi}, \hat{p})_{H^0} - (v_\varepsilon, \widehat{A}^*\hat{p})_{H^0} + (A_\varepsilon v_\varepsilon + \varphi_\varepsilon, \hat{p})_{H^0}\} dt =$$

$$= \int_0^T \{(v_\varepsilon, A_\varepsilon^*\hat{p})_{H^0} - (v_\varepsilon, \widehat{A}^*\hat{p})_{H^0} + (\varphi_\varepsilon - \hat{\varphi}, \hat{p})_{H^0}\} dt. \quad (2.8)$$

The left and right sides of (2.8) are clearly also equal for $\varphi_\varepsilon \in X_1^0$ and $\hat{\varphi} \in X_1^0$, because in this case the two sides contain $A_\varepsilon v_\varepsilon$, $\widehat{A}\hat{v}$, etc., in which there are derivatives of $v_\varepsilon$ or $\hat{v}$, and $\varphi_\varepsilon$ and $\hat{\varphi}$ can be approximated in $X_1^0$ by functions from $X^1$.

The relation

$$\int_0^T \{(v_\varepsilon, A_\varepsilon^*\hat{p})_{H^0} - (v_\varepsilon, \widehat{A}^*\hat{p})_{H^0}\} dt = \int_0^T \{(\hat{v}, A_\varepsilon^*\hat{p})_{H^0} - (\hat{v}, \widehat{A}^*\hat{p})_{H^0} + (v_\varepsilon - \hat{v}, (A_\varepsilon^* - \widehat{A}^*)\hat{p})_{H^0}\} dt \quad (2.9)$$

holds, and (2.8) and (2.9) imply that

$$\varepsilon^{-1}\{\Phi[u(\cdot, \varepsilon|\mu)] - \Phi[\hat{u}(\cdot)]\} = \varepsilon^{-1} \int_0^T \{(\hat{v}(x, t), [A_\varepsilon^*(x, t) -$$

$$- \widehat{A}(x, t)]\hat{p}(x, t))_{H^0} + (\varphi_\varepsilon(x, t) - \hat{\varphi}(x, t), \hat{p}(x, t))_{H^0} + \xi(\varepsilon)\} dt = \varepsilon^{-1} \int_{\Theta(\varepsilon)} \eta(t, \bar{u}) dt + \xi(\varepsilon); \quad (2.10)$$

here $\eta(t, \bar{u}) = (\hat{v}(x, t), [A^*(x, \bar{u}, t) - A^*(x, \hat{u}(t), t)]\hat{p}(x, t))_{H^0} + (\varphi(x, \bar{u}, t) - \varphi(x, \hat{u}(t), t), \hat{p}(x, t))_{H^0}$, and

$$\xi(\varepsilon) = \xi(\varepsilon, \mu) = \varepsilon^{-1} \int_0^T (v_\varepsilon(x, t) - \hat{v}(x, t), [A_\varepsilon^*(x, t) - \widehat{A}(x, t)]\hat{p}(x, t))_{H^0} dt. \quad (2.11)$$

LEMMA 2.1. The quantity (2.11) tends to zero when $\varepsilon \to +0$, for $\mu = \|\bar{t}, \bar{u}\| \in (0, T) \times \Delta$.

Suppose that this lemma has been proved. Then the quantity (2.10) has a limit when $\varepsilon \to +0$ for all Lebesgue points $\bar{t} \in (0, T)$ if $\eta(t, \bar{u}): [0, T] \to \mathbb{R}$, i.e., for $\bar{t} \in \Theta_1 \cap \Theta_2(\bar{u}) \cap \Theta_3(\bar{u})$, where $\Theta_3(\bar{u})$ is the Lebesgue-point set of $q(t) = (\hat{v}(x, t), A^*(x, \bar{u}, t)\hat{p}(x, t))_{H^0}$. Theorem 4.2 [3] implies that $\hat{p}(x, t)$ is a weakly continuous function of $t$ in $H^1$; hence, by virtue of our assumptions concerning $f$ and $\partial f/\partial x$, we have $\Theta_3(\bar{u}) = [0, T]$ ($\forall \bar{u}$). The expression in (2.10) has a limit for $\varepsilon \to +0$ when $\bar{t} \in \Theta_1 \cap \Theta_2(\bar{u})$, and is nonnegative because $\hat{u}(\cdot)$ is optimum (this limit is the required differential of the functional (2.3) on the curve $u(\cdot, \varepsilon|\mu)$ [3]). The conclusion of the theorem therefore follows.

Proof of Lemma 2.1. It follows from the construction of the bundle of variations $u(\cdot, \varepsilon|\mu)$ that

$$\xi(\varepsilon) = \varepsilon^{-1} \int_{\Theta(\varepsilon)} (v_\varepsilon(x, t) - \hat{v}(x, t), \zeta(x, t))_{H^0} dt,$$

where
$$\zeta(x, t) = \zeta(x, t, \mu) = (A^*(x, \bar{u}, t) - \hat{A}^*(x, t))\hat{p}(x, t), \qquad (2.12)$$
and $\zeta(x, t) \in X^0$. Theorem 4.2.1 in [3] implies that $\hat{p} \in L^\infty([0, T] \to H^1)$ and $\zeta \in L^\infty([0, T] \to H^0)$.

Suppose that $q_\varepsilon(\cdot, \cdot, \rho) = \mathscr{L}_\rho^*[u(\cdot, \varepsilon|\mu)]\zeta(\cdot, \rho)$ and $\hat{q}(\cdot, \cdot, \rho) = \mathscr{L}_\rho^*[\hat{u}(\cdot)]\zeta(\cdot, \rho)$; these functions are defined for $\rho$ such that $\zeta(\cdot, \rho) \in H^0$, i.e., for almost all $\rho \in [0, T]$. Lemma 1.2 implies that

$$\zeta(\varepsilon) = \varepsilon^{-1} \int_{\Theta(\varepsilon)} d\rho \int_\rho^T \{(\varphi_\varepsilon(x, t), q_\varepsilon(x, t, \rho))_{H^0} - (\hat{\varphi}(x, t), \hat{q}(x, t, \rho))_{H^0}\} dt = \gamma_1(\varepsilon) + \gamma_2(\varepsilon), \qquad (2.13)$$

where
$$\gamma_1(\varepsilon) = \gamma_1(\varepsilon, \mu) = \varepsilon^{-1} \int_{\Theta(\varepsilon)} d\rho \int_\rho^{\bar{t}+\varepsilon} (\varphi(x, \bar{u}, t) - \hat{\varphi}(x, t), q_\varepsilon(x, t))_{H^0} dt,$$

and
$$\gamma_2(\varepsilon) = \gamma_2(\varepsilon, \mu) = \varepsilon^{-1} \int_{\Theta(\varepsilon)} d\rho \int_\rho^T (\hat{\varphi}(x, t), q_\varepsilon(x, t, \rho) - \hat{q}(x, t, \rho))_{H^0} dt.$$

It follows from Lemma 1.1 that there is a constant $c > 0$ for which $\|q_\varepsilon(x, t, \rho)\|_{\mathscr{C}^0[\rho, T]} \leq c\|\zeta(x, \rho)\|_{H^0}$ ($\forall \varepsilon \in [0, \varepsilon_\mu]$). Hence

$$|\gamma_1(\varepsilon)| \leq c\varepsilon^{-1} \int_{\bar{t}}^{\bar{t}+\varepsilon} \{\|\zeta(x, \rho)\|_{H^0} \int_\rho^{\bar{t}+\varepsilon} \|\varphi(x, \bar{u}, t) - \hat{\varphi}(x, t)\|_{H^0} dt\} d\rho \leq$$
$$\leq (c \int_{\bar{t}}^{\bar{t}+\varepsilon} \|\varphi(x, \bar{u}, t) - \hat{\varphi}(x, t)\|_{H^0} dt) \varepsilon^{-1} \int_{\bar{t}}^{\bar{t}+\varepsilon} \|\zeta(x, \rho)\|_{H^0} d\rho \to 0 \qquad (2.14)$$

by virtue of condition II and the essential boundedness of $\|\zeta(x, \rho)\|_{H^0}$ as a function of $\rho$.

It remains to prove that $\gamma_2(\varepsilon) = \gamma_2(\varepsilon, \mu) \to 0$ when $\varepsilon \to 0$ for arbitrary $\mu$. Let $w_\varepsilon = L[u(\cdot, \varepsilon|\mu)]\hat{\varphi} - L[\hat{u}(\cdot)]\hat{\varphi}$; then Lemma 1.2 implies that

$$\gamma_2(\varepsilon) = \varepsilon^{-1} \int_{\Theta(\varepsilon)} (w_\varepsilon(x, t), \zeta(x, t))_{H^0} dt,$$

and
$$|\gamma_2(\varepsilon)| \leq \operatorname*{ess\,sup}_{t \in \Theta(\varepsilon)} (w_\varepsilon(x, t), \zeta(x, t))_{H^0} \varepsilon^{-1} \int_{\Theta(\varepsilon)} dt. \qquad (2.15)$$

Lemma 2.1 will now have been proved if we show that

$$\operatorname*{ess\,sup}_{t \in [0, T]} |(w_\varepsilon(x, t), \zeta(x, t))_{H^0}| \to 0 \text{ for } \varepsilon \to +0 \ (\forall \mu). \qquad (2.16)$$

**Proposition 2.1.** If $\bar{\zeta}(x, \rho) \in L^\infty([0, T] \to H^0)$, and there is a bounded domain $D_0 \subset \mathbf{R}^n$ such that, for almost all $\rho \in [0, T]$ we have $\bar{\zeta}(x, \rho) = 0$ for almost all $x \in \mathbf{R}^n \setminus D_0$. Then $\operatorname*{ess\,sup}_{\rho \in [0, T]} |(w_\varepsilon(x, \rho), \bar{\zeta}(x, \rho))_{H^0}| \to 0$ for $\varepsilon \to +0$.

**Proof.** We need to prove this only for almost-everywhere nonnegative functions $\bar{\zeta}(x, \rho)$, and then apply this result in the general case to nonnegative $\bar{\zeta}^+$ and $\bar{\zeta}^-$ such that $\bar{\zeta}(x, \rho) = \bar{\zeta}^+(x, \rho) - \bar{\zeta}^-(x, \rho)$. We therefore assume that $\bar{\zeta}(x, \rho) \geq 0$. We also assume that $\int_{\mathbf{R}^n} \bar{\zeta}(x, \rho) dx = 1$ for almost all $\rho$ (the transition to the general case is obvious).

Let $\mathfrak{T} \subset [0, T]$ be the set of $\rho \in [0, T]$ for which $\bar{\zeta}(x, \rho) \in H^0$, $\int_{\mathbf{R}^n} \bar{\zeta}(x, \rho) dx = 1$, $\bar{\zeta}(x, \rho) = 0$ for $x \in \mathbf{R}^n \setminus D$, and $\bar{\zeta}(x, \rho) \geq 0$ otherwise. Let $b(\rho)$ be a random variable given for $\rho \in \mathfrak{T}$, such that its density distribution is $\bar{\zeta}(x, \rho)$. Then $M|b(\rho)| \in L_\infty(0, T)$ and mes $\mathfrak{T} = $ mes $[0, T]$. Let $y_\varepsilon^{b, \rho}(t)$ and $\hat{y}^{b, \rho}(t)$ be the solutions of (2.1) on $[\rho, T]$, with $u(\cdot) = u(\cdot, \varepsilon|\mu)$ and $u(\cdot) = \hat{u}(\cdot)$, respectively, and with the initial condition $y(\rho) = b(\rho)$. Theorem 1.2 implies that $q_\varepsilon(x, t, \rho)$ and $\hat{q}(x, t, \rho)$ are the distribution densities of $y_\varepsilon^{b, \rho}(t)$ and $\hat{y}^{b, \rho}(t)$, respectively.

There is clearly a sequence of bounded continuous functions $\varphi_i(x, t)$, satisfying a global Lipschitz condition and such that $\varphi_i \to \hat{\varphi}$ for $i \to +\infty$ in the metric of $X_1^0$. Let $w_\varepsilon^i = L[u(\cdot, \varepsilon|\mu)]\varphi^{(i)} - L[\hat{u}(\cdot)]\varphi^{(i)}$. By virtue of the customary reasoning used in connection with Gronwall's lemma, we conclude, for $\rho \in \mathfrak{T}$, that

$$|(w_\varepsilon^i(x, \rho), \bar{\zeta}(x, \rho))_{H^0}| = |\int_\rho^T (q_\varepsilon(x, t, \rho) - \hat{q}(x, t, \rho), \varphi^{(i)}(x, t))_{H^0} dt| \leq$$

$$\leqslant |M \int_\rho^T \{\varphi^{(i)}[y_\varepsilon^{b,\rho}(t), t] - \varphi^{(i)}[\hat{y}^{b,\rho}(t), t]\} dt| \leqslant K_{\varphi^{(i)}} M(|b(\rho)|+1)\varepsilon. \quad (2.17)$$

Here the $K_{\varphi^{(i)}} > 0$ are constants independent of $\varphi^{(i)}$.

It follows from Lemma 1.1 that
$$\sup_{\rho \in [0, T]} \|w_\varepsilon^i(x, \rho) - w_\varepsilon(x, \rho)\|_{H^0} \leqslant 2c_\mu \|\varphi^{(i)} - \hat{\varphi}\|_{X_1^0}. \quad (2.18)$$

where $c_\mu$ is the supremum, over $\varepsilon \in [0, \varepsilon_\mu]$, of the norm of the operator $L[u(\cdot, \varepsilon | \mu)]: X_1^0 \to \mathscr{C}^0$. Lemma 1.1 implies that $\sup_\mu c_\mu < +\infty$.

Let $\delta > 0$ be arbitrary. We prove that there is $\varepsilon_* > 0$ such that $\operatorname*{ess\,sup}_{\rho \in [0, T]} |(w_\varepsilon(x, \rho), \bar{\zeta}(x, \rho))_{H^0}| < \delta$ ($\forall \varepsilon \in [0, \varepsilon_\mu]$), and this will prove Proposition 2.1.

Let $\alpha_\zeta = \operatorname*{ess\,sup}_{\rho \in [0, T]} \|\bar{\zeta}(x, \rho)\|_{H^0}$. When $\alpha_\zeta = 0$ the conclusion of Proposition 2.1 is obvious. If $\alpha_\zeta > 0$, then there is clearly an index $i_*$ such that $\|\varphi^{(i)} - \hat{\varphi}\|_{X_1^0} \leqslant \delta/4c_\mu \alpha_\zeta$ for $i \geqslant i_*$.

The number $\varepsilon_*$ is selected so that
$$(w_\varepsilon^{i_*}(x, \rho), \bar{\zeta}(x, \rho))_{H^0} < \delta/2 \quad (\forall \rho \in \mathfrak{T}, \forall \varepsilon \in [0, \varepsilon_*]), \quad (2.19)$$

which is possible by virtue of (2.17). From (2.18) and (2.19) we obtain
$$|(w_\varepsilon(x, \rho), \bar{\zeta}(x, \rho))_{H^0}| \leqslant |(w_\varepsilon(x, \rho) - w_\varepsilon^{i_*}(x, \rho), \bar{\zeta}(x, \rho))_{H^0}| +$$
$$+ |(w_\varepsilon^{i_*}(x, \rho), \bar{\zeta}(x, \rho))_{H^0}| \leqslant \frac{\delta}{4c_\mu \alpha_\zeta} 2c_\mu \alpha_\zeta + \frac{\delta}{2} = \delta \quad (\forall \rho \in \mathfrak{T}).$$

This proves Proposition 2.1 because $\operatorname{mes} \mathfrak{T} = \operatorname{mes}[0, T]$.

We continue the proof of Lemma 2.1. It remains to show that (2.16) holds, and then the required result follows from (2.13)-(2.15).

Let $D_i \subset \mathbf{R}^n$ be a sequence of bounded domains, and let the functions $\rho_i \in H^1$, $i = 1, 2, \ldots$ be such that $\rho_i(x) = 0$ for almost all $x \in \mathbf{R}^n \setminus D_i$, and $\rho_i(x) \to \rho_0(x)$ in the $H^1$-metric for $i \to +\infty$ (it is easily shown that such a sequence exists). Let $p_i = \mathscr{L}_0^*[\hat{u}(\cdot)]\rho_i$ and let $\bar{\zeta}_i(x, t) = \bar{\zeta}_i(x, t, \mu) = (A^*(x, \bar{u}, t) - \hat{A}^*(x, t)) p_i(x, t)$ [see (2.12)]. It follows from Theorem 4.2.1 in [3] that $p_i \in L^\infty([0, T] \to H^1)$ and $p_i \to p$ in the metric of $L^\infty([0, T] \to H^1)$, and $\bar{\zeta}_i \in L^\infty([0, T] \to H^0)$ and $\bar{\zeta}_i \to \zeta$ in the metric of $L^\infty([0, T] \to H^0)$ for $i \to +\infty$, uniformly with respect to $\varepsilon$ and $\mu$. Moreover $\bar{\zeta}_i(x, t) = 0$ for almost all $x \in \mathbf{R}^n \setminus \tilde{D}_i$, where the $\tilde{D}_i \subset \mathbf{R}^n$ are bounded domains; this is clear from Theorem 1.2 and the boundedness of f in (2.1).

Let $\delta > 0$ be arbitrary, and assume that $\|\hat{\varphi}\|_{X_1^0} > 0$ [otherwise (2.16) is not needed in the proof].

There is clearly an index $i_*$ for which
$$\operatorname*{ess\,sup}_{\rho \in [0, T]} \|\bar{\zeta}_{i_*}(x, \rho) - \zeta(x, \rho)\|_{H^0} \leqslant (2c_\mu \|\hat{\varphi}\|_{X_1^0})^{-1} \frac{\delta}{2} \quad (\forall \varepsilon \in [0, \varepsilon_\mu]).$$

Here $c_\mu$ is the same constant as in (2.18), and is independent of $\varepsilon$. There is $\varepsilon_* > 0$ such that
$$\operatorname*{ess\,sup}_{\rho \in [0, T]} |(w_\varepsilon(x, \rho), \bar{\zeta}_{i_*}(x, \rho))_{H^0}| \leqslant \frac{\delta}{2} \quad (\forall \varepsilon \in [0, \varepsilon_*])$$

by virtue of Proposition 2.1. Hence Lemma 1.1 implies that, for $\varepsilon \in [0, \varepsilon_*]$, the following relations hold:
$$\operatorname*{ess\,sup}_{\rho \in [0, T]} |(w_\varepsilon(x, \rho), \zeta(x, \rho))_{H^0}| \leqslant \operatorname*{ess\,sup}_{\rho \in [0, T]} |(w_\varepsilon(x, \rho), \bar{\zeta}_{i_*}(x, \rho))_{H^0}| +$$
$$+ \operatorname*{ess\,sup}_{\rho \in [0, T]} |(w_\varepsilon(x, \rho), \zeta(x, \rho) - \bar{\zeta}_{i_*}(x, \rho))_{H^0}| \leqslant \frac{\delta}{2} + \frac{\delta}{2} = \delta.$$

This proves (2.16) and completes the proof of Lemma 2.1 and Theorem 2.1.

We next give an existence theorem for an optimal solution of a nonsmooth optimal-control problem.

Let $\Delta$ be a compact in $\mathbf{R}^m$, let $C^*(\Delta)$ be the conjugate to the Banach space $C(\Delta \to \mathbf{R})$ of continuous functions on $\Delta$, and assume that $C^*(\Delta)$ is endowed with the weak norm. We write rpm($\Delta$) for the set of regular probability measures on $\Delta$ ([9], p. 65); then $\text{rpm}(\Delta) \subset C^*(\Delta)$. We introduce the class $\mathcal{U}_R$ of admissible controls, containing functions $u \in L^\infty([0, T] \to C^*(\Delta))$ such that $u(t) \in \text{rpm}(\Delta)$ for almost all $t \in [0, T]$. Suppose that $f(x, u, t): \mathbf{R}^n \times \text{rpm}(\Delta) \times [0, T] \to \mathbf{R}^n$ in (2.1) and $f(x, u, t) = \int_\Delta f(x, v, t) u(dv)$ for $u \in \text{rpm}(\Delta)$. For $u(\cdot) \in \mathcal{U}_R$ a solution of the Cauchy problem (2.1), (2.2) is unique, and for nonrandom a the solution is nonrandom.

THEOREM 2.2. Suppose that the random vector a has distribution density $\rho_0(x) \in L_2(\mathbf{R}^n)$, the function $f: \mathbf{R}^n \times \mathbf{R}^m \times [0, T] \to \mathbf{R}^n$ satisfies the conditions of Theorem 2.1, $\varphi(x, u, t) \equiv \varphi(x, t) \in X_1^0$ is independent of u, and the set $\Delta$ is compact in $\mathbf{R}^m$. Then problem (2.1)-(2.3) has an optimal control $\hat{u}(\cdot)$ in $\mathcal{U}_R$.

COROLLARY 2.1. If $f(x, u, t)$ depends linearly on $u \in \Delta$ for all x and t, and $\Delta$ is convex, then under the assumptions of Theorem 2.2 there is an optimal control in $\mathcal{U}_0$ (i.e., an ordinary control).

Proof of Theorem 2.2. It is known ([9], p. 304) that the set $\mathcal{U}_R$ is convex, compact, and sequentially compact if it is endowed with the following topology: for a sequence $u_j(\cdot) \in \mathcal{U}_R$ and $\tilde{u}(\cdot) \in \mathcal{U}_R$, we have $u_i(\cdot) \to \tilde{u}(\cdot)$ for $i \to +\infty$ if

$$\int_0^T dt \int_\Delta \alpha(t, v) u_i(t, dv) \to \int_0^T dt \int_\Delta \alpha(t, v) \tilde{u}(t, dv)$$

when $\alpha(t, v) \in L^1([0, T] \to C(\Delta \to \mathbf{R}))$. To prove the theorem it is sufficient to show that, for such a sequence, $\Phi[u_i(\cdot)] \to \Phi[\tilde{u}(\cdot)]$ when $i \to +\infty$ [i.e., the functional (2.3) is sequentially continuous in the topology introduced].

Let $y_i(\cdot)$ be the solutions of the Cauchy problem (2.1), (2.2) corresponding to $u(\cdot) = u_i(\cdot)$, and let $\tilde{y}(\cdot)$ be its solution corresponding to $u(\cdot) = \tilde{u}(\cdot)$, and let $p_i = \mathscr{L}_0^*[u_i(\cdot)]\rho_0$ and $\tilde{p} = \mathscr{L}_0^*[\tilde{u}(\cdot)]\rho_0$. It is known ([9], p. 358) that $y_i(t) \to \tilde{y}(t)$ with probability 1, uniformly with respect to $t \in [0, T]$ for $i \to +\infty$. Theorem 1.2 implies that $p_i(x, t)$ and $\tilde{p}(x, t)$ are the distribution densities of $y_i(t)$ and $\tilde{y}(t)$, respectively; hence $p_i \to \tilde{p}$ weakly in $X^0$ for $i \to +\infty$. Moreover $(p_i, \varphi)_{X^0} \to (\tilde{p}, \varphi)_{X^0}$ for $i \to +\infty$, $\Phi[u_i(\cdot)] = (p_i, \varphi)_{X^0}$, and $\Phi[\tilde{u}(\cdot)] = (\tilde{p}, \varphi)_{X^0}$. Hence $\Phi[u_i(\cdot)] \to \Phi[\tilde{u}(\cdot)]$ for $i \to +\infty$, and this proves Theorem 2.2.

It is natural to extend our necessary conditions for optimization to the admissible-control class $\mathcal{U}_R$.

THEOREM 2.3. If $\hat{u}(\cdot) \in \mathcal{U}_R$ and $\rho_0(x) \in W_2^1(\mathbf{R}^n)$, then all conclusions of Theorem 2.1 remain in force: if $\bar{u} \in \Delta$ then condition (2.5) is satisfied for almost all t; if (2.6) is satisfied for almost all t, then (2.5) holds for $\bar{u} \in \Delta$.

In Theorem 2.3, (2.5) holds for $\bar{u} \in \text{rpm}(\Delta)$. However this is superfluous because it does not strengthen the conclusion [the upper bounds with respect to $\bar{u} \in \Delta$ and $\bar{u} \in \text{rpm}(\Delta)$ coincide].

The proof of Theorem 2.3 is the same as the proof of Theorem 2.1 with $\mathcal{U}_0$ replaced by $\mathcal{U}_R$.

We conclude by noting that a similar "regularization" of a central problem by restrictions on the summability of the initial density is used in [10], for a special control problem for nonsingular stochastic Ito equations. The results of the present article were announced in [11].

LITERATURE CITED

1. R. Courant, Partial Differential Equations [Russian translation], Moscow (1964).
2. V. F. Dem'yanov, T. K. Vinogradova, V. N. Nikulina, et al., Nonsmooth Problems in Optimization and Control Theory [in Russian], Leningrad (1982).
3. B. L. Rozovskii, Evolution Stochastic Systems. Linear Theory with Applications to Static Random Processes [in Russian], Moscow (1983).
4. O. A. Ladyzhenskaya, Boundary-Value Problems of Mathematical Physics [in Russian], Moscow (1973).
5. A. N. Kolmogorov and S. V. Fomin, Elements of the Theory of Functions and Functional Analysis [in Russian], Moscow (1981).
6. N. Bourbaki, Integration (Measures, Integration of Measures) [Russian translation], Moscow (1967).
7. N. G. Dokuchaev and V. A. Yakubovich, Cybernetics and Numerical Techniques [in Russian], No. 54, Kiev (1983), pp. 72-78.

8. V. A. Yakubovich, Sib. Mat. Zh., 18, No. 3, 685-707 (1977).
9. G. Varga, Optimum Control for Differential and Functional Equations [Russian translation], Moscow (1972).
10. N. G. Dokuchaev, Differents. Uravn., 27, No. 3, 403-414 (1991).
11. N. G. Dokuchaev, Proceedings of the 3rd All-Union School, "Pontryagin Lectures. Optimum Control. Geometry and Analysis," Kemerovo (1990), p. 127.

# CONTROLLABILITY SET OF AN ALMOST-PERIODIC LINEAR SYSTEM

A. G. Ivanov and E. L. Tonkov

UDC 517.977.1

**1. Definitions and Notation.** We use the following notation: $R^n$ denotes n-dimensional Euclidean space; $\langle x, y \rangle$ is the scalar product of $x, y \in R^n$ and $|x| = \sqrt{\langle x, x \rangle}$, $\text{Hom}(R^n)$ is the space of linear operators $A: R^n \to R^n$ with the norm $|A| \doteq \max_{|x|=1} |Ax|$; $\text{conv}(R^n)$ is the space of compact convex sets in $R^n$ with the Hausdorff metric $\text{dist}(\cdot, \cdot)$, $h(\psi, F) \doteq \max_{f \in F} \langle \psi, f \rangle$ is the support function of a set $F \in \text{conv}(R^n)$ in the direction of the vector $\psi$. Properties of support functions, established for example in [1], will be used without referring to their existence.

We write $(\mathfrak{A}_l, \rho_1)$ $(l > 0)$ for the metric space of locally summable mappings of $R$ into $\text{Hom} \times (R^n)$, with the metric $\rho_1(A, B) < \infty$ defined by the formula

$$\rho_1(A, B) \doteq \sup_t \frac{1}{l} \int_t^{t+l} |A(s) - B(s)| ds, \quad t \in R, \; A, B \in \mathfrak{A}_l,$$

and $(\mathfrak{B}_l, \rho_2)$ for the metric space of mappings $V: R \to \text{conv}(R^n)$ such that $t \mapsto |V(t)| \doteq \text{dist}(\{0\}, V(t))$ is measurable (all measures are Lebesgue measures), $\operatorname{ess\,sup}_{t \in R} |V(t)| \doteq k_V < \infty$, and

$$\rho_2(V, W) \doteq \sup_t \frac{1}{l} \int_t^{t+l} \text{dist}(V(s), W(s)) ds, \quad t \in R, \; V, W \in \mathfrak{B}_l.$$

We set up a correspondence between each pair $(A, V) \in \mathfrak{A}_l \times \mathfrak{B}_l$ and the function $t \mapsto \xi(t) \doteq (A(t), V(t))$, and on $\mathfrak{A}_l \times \mathfrak{B}_l$ we define the metric $d_l(\xi, \eta) \doteq \rho_1(A, B) + \rho_2(V, W)$, where $\xi = (A, V)$ and $\eta = (B, W)$.

In the following $f_\tau(\cdot) \doteq f(\cdot + \tau)$ denotes the translation, through the distance $\tau \in R$ of $f$, and $f_0 \doteq f$.

A function $\xi \in \mathfrak{A}_l \times \mathfrak{B}_l$ is said to be Stepanov almost periodic (a.p.) if, for some $\varepsilon > 0$, the set $\mathscr{E}_l(\xi, \varepsilon) \doteq \{\tau \in R : d_l(\xi_\tau, \xi) < \varepsilon\}$ of its $\varepsilon$-almost periods ($\varepsilon$-a.p.) is relatively everywhere dense in $R$. The set of Stepanov a.p. functions from $\mathfrak{A}_l \times \mathfrak{B}_l$ will be denoted by $S_l(R, \mathscr{H})$, where $\mathscr{H} \doteq \text{Hom}(R^n) \times \text{conv}(R^n)$. It can be proved [2] that, if $\xi \in S_l(R, \mathscr{H})$, then it is $d_\ell$-bounded [i.e., $d_\ell(\xi, 0) < \infty$] and $d_\ell$-uniformly continuous [i.e., if $\varepsilon > 0$ then there is $\delta > 0$ such that $d_\ell(\xi_\tau, \xi) < \varepsilon$ for $|\tau| < \delta$].

It can be directly verified that $d_\ell$-distances for distinct $\ell > 0$ are topologically equivalent, since

$$d_{l_1}(\xi, \eta) \geq 2 d_{l_2}(\xi, \eta), \quad d_{l_2}(\xi, \eta) \geq \frac{l_1}{l_2} d_{l_1}(\xi, \eta). \tag{1.1}$$

for $\ell_1 < \ell_2$. We therefore put $\ell = 1$ and write simply d-distance, $\mathfrak{A}, \mathfrak{B}$, etc., omitting the subscript.

We also recall that a continuous function f on R into a metric space $(\mathscr{X}, \rho)$, is called Bohr a.p. if, for each $\varepsilon > 0$, the set $E(f, \varepsilon) \doteq \{\tau \in R : \sup_{t \in R} \rho(f_\tau(t), f(t)) < \varepsilon\}$ of its $\varepsilon$-a.p. is relatively everywhere dense on R.